
\magnification 1250
\font\san=cmssdc10
\def\sy{\hbox{\san\char 83}} 

\mathcode`A="7041 \mathcode`B="7042 \mathcode`C="7043
\mathcode`D="7044 \mathcode`E="7045 \mathcode`F="7046
\mathcode`G="7047 \mathcode`H="7048 \mathcode`I="7049
\mathcode`J="704A \mathcode`K="704B \mathcode`L="704C
\mathcode`M="704D \mathcode`N="704E \mathcode`O="704F
\mathcode`P="7050 \mathcode`Q="7051 \mathcode`R="7052
\mathcode`S="7053 \mathcode`T="7054 \mathcode`U="7055
\mathcode`V="7056 \mathcode`W="7057 \mathcode`X="7058
\mathcode`Y="7059 \mathcode`Z="705A
\font\eightrm=cmr8         \font\eighti=cmmi8
\font\eightsy=cmsy8        \font\eightbf=cmbx8
        
        \font\sixrm=cmr6
\font\sixi=cmmi6           \font\sixsy=cmsy6
\font\sixbf=cmbx6
\def\eightpoint{%
  \textfont0=\eightrm \scriptfont0=\sixrm \scriptscriptfont0=\fiverm
  \def\rm{\fam0\eightrm}%
  \textfont1=\eighti  \scriptfont1=\sixi  \scriptscriptfont1=\fivei
  \textfont2=\eightsy \scriptfont2=\sixsy \scriptscriptfont2=\fivesy
  \textfont\bffam=\eightbf \scriptfont\bffam=\sixbf
  \scriptscriptfont\bffam=\fivebf
  \def\bf{\fam\bffam\eightbf}%
   \normalbaselineskip=9pt\rm}
\def\spacedmath#1{\def\packedmath##1${\bgroup\mathsurround
=0pt##1\egroup$} \mathsurround#1
\everymath={\packedmath}\everydisplay={\mathsurround=0pt}}
\def\nospacedmath{\mathsurround=0pt
\everymath={}\everydisplay={} } \spacedmath{2pt}
\def\phfl#1#2{\normalbaselines{\baselineskip=0pt
\lineskip=10truept\lineskiplimit=1truept}\nospacedmath\smash
{\mathop{\hbox to 8truemm{\rightarrowfill}}
\limits^{\scriptstyle#1}_{\scriptstyle#2}}}
\def\diagram#1{\def\normalbaselines{\baselineskip=0truept
\lineskip=10truept\lineskiplimit=1truept}   \matrix{#1}}
\def\vfl#1#2{\llap{$\scriptstyle#1$}\left\downarrow\vbox
to 6truemm{}\right.\rlap{$\scriptstyle#2$}}
\def\eqalign#1{\null\,\vcenter{\openup\jot\ialign{
\strut\hfil$\displaystyle{##}$&$\displaystyle{{}##}$\hfil
&&\quad\strut\hfil$\displaystyle{##}$&$\displaystyle{{}##}$\hfil
\crcr#1\crcr}}\,}

\def\up#1{\raise 1ex\hbox{\sevenrm#1}}
\def\tx{\kern-1.5pt -}
\def\cqfd{\kern 2truemm\unskip\penalty 500\vrule height
4pt depth 0pt width 4pt\medbreak}  \def\virg{\raise .4ex\hbox{,}}
\def\decale#1{\smallbreak\hskip 28pt\llap{#1}\kern 5pt}
\def\no{n\up{o}\kern 2pt}
\def\ind{\par\hskip 1truecm\relax}
\def\indp{\par\hskip 0.5truecm\relax}
\parindent=0cm

\def\rond{\kern 1pt{\scriptstyle\circ}\kern 1pt}

\vsize = 25truecm
\hsize = 16truecm
\hoffset = -.15truecm
\voffset = -.5truecm
\baselineskip15pt
\overfullrule=0pt
\centerline{\bf Quantum cohomology of complete intersections}
\medskip
\centerline{Arnaud Beauville{\footnote{\parindent
0.4cm$^1$}{\vtop{\hsize15.5truecm\noindent \eightpoint
Partially supported by the
European HCM project ``Algebraic Geometry in Europe" (AGE).
}}\parindent 0cm}}
\vskip1.1truecm
{\bf Introduction}
\smallskip
\ind The quantum cohomology algebra of a projective manifold $X$ is the
coho\-mo\-logy of $X$ endowed with a different algebra structure,
which takes into account the geometry of rational curves in $X$.  This
structure has been first defined heuristically by the mathematical physicists
[V,W]; a rigorous construction (and proof of the associativity, which is highly
non trivial) has been achieved recently by Ruan and Tian [R-T].

\ind When computed e.g. for surfaces,  the quantum cohomology looks
rather complicated [C-M]. The aim of this note is to show that the situation
improves consi\-derably when the dimension becomes high with respect to
the degree. Our main result is:
\smallskip
{\bf Theorem}$.-$  {\it Let $X\i {\bf
P}^{n+r}$ be a smooth complete intersection of degree $(d_1,\ldots,d_r)$ and
dimension
$n\ge 2$, with
 $n\ge 2\sum (d_i-1)-1$. The quantum cohomology algebra $H^*(X,{\bf Q})$ is
the algebra generated by the hyperplane class $H$ and the primitive
cohomology $H^n(X,{\bf Q})_{\rm o}$ , with the  relations:
$$H^{n+1}=d_1^{d_1}\ldots d_r^{d_r}\,H^{d-1}\qquad H\cdot\alpha=0\qquad
\alpha\cdot\beta=(\alpha\,|\,\beta){1\over d} \,(H^n-d_1^{d_1}\ldots
d_r^{d_r}\,H^{d-2})$$ for $\alpha,\beta\in H^n(X,{\bf Q})_{\rm o}$} .

\ind The
method applies more generally to a large class of Fano manifolds (see
Proposition 1
below). It is actually a straightforward consequence of the definitions
-- except for the exact value of the coefficient $d_1^{d_1}\ldots
d_r^{d_r}$, which requires some standard
computations in the cohomology of the Grassmannian. Still I believe
that the simplicity of the result is worth noticing.
As usual it implies a number of enumerative formulas:
for instance we find that the number of conics passing through $2$
general points in a hypersurface of degree $d$ and dimension $2d-3$ is
${1\over 2}d!(d-1)!$, while the number of twisted cubics through $3$
general points in a
hypersurface of  degree $d$ and dimension $3d-6$ is $d!((d-1)!)^2$.
\vskip1truecm
{\bf 1. Quantum cohomology of Fano manifolds}
\smallskip
\ind I am considering in this paper Fano manifolds  with
$b_2=1$, i.e. smooth compact complex manifolds $X$ such that $H^2(X,{\bf
Z})$ is generated by an ample class $H$ and the canonical class $K_X$ is
$-kH$ for some positive integer $k$.
I will use the follo\-wing properties of the  quantum cohomo\-logy product   on
$H^* (X,{\bf Z})$ (proved in [R-T]):

(1.1) \  it is invariant under smooth deformations;

(1.2) \  it is associative, compatible with the grading mod.$2$, and
anticommutative. It has the same unit $1\in H^0(X,{\bf Z})$ and the same
intersection form $(\ |\ )$ as the ordinary cohomology;

(1.3) \  the product $x\cdot y$  of two homogeneous elements is
defined by $$x\cdot y = (x\cdot y)_0+(x\cdot y)_1+\ldots+(x\cdot
y)_j+\ldots$$ where $(x\cdot y)_0$ is the ordinary cohomology
product, and $(x\cdot y)_j$ is a class of degree $\deg(x)+\deg(y)-2kj$.

(1.4) \  Assume that the  moduli space ${\cal
M}_j$ of maps $f:{\bf P}^1\rightarrow X$ of
degree $j$ (i.e. such that $\deg f^*H=j$) has the expected dimension $n+kj$;
choose any smooth compactification $\overline{{\cal M}}_j$ of ${\cal M}_j$
such that the evaluation maps  $e_i:{\cal M}_d\rightarrow X$ ($0\le i\le 2$)
defined by $e_i(f)=f(i)$   extend to $\overline{{\cal M}}_j$.
Then the ``instanton correction" $(x\cdot y)_j$ is
defined by $$\langle x,y,z\rangle^{}_j:=\bigl((x\cdot y)_j\ |\
z\bigr)=\int_{\overline{\cal M}_j}e_0^*x\,.\,e_1^*y\,.\,e_2^*z \ .$$
(1.5) \  If $x,y,z\in H^*(X,{\bf Z})$
are classes of  subvarieties $A,B,C$ of $X$ which are in general position, it
follows easily from (1.4) that the triple product  $\langle x,y,z\rangle^{}_j$
is  the number of  curves of degree $j$ meeting $A$, $B$ and $C$ (counted with
multiplicity $abc$ if the curve meets $A$, resp. $B$, resp. $C$ in $a$,
resp. $b$, resp. $c$ distinct points). \smallskip

\ind To avoid confusion I will denote by $H_p\in H^{2p}(X,{\bf Z})$  for $0\le
p\le n$ the $p$\tx th power of $H$ in the ordinary cohomology, and reserve  the
notations $x\cdot y$ or $x^n$ ($x,y\in H^*(X,{\bf Q})$) exclusively for the
quantum product. One has $H_0=1$,
$H_1=H$, and $H_n$ is $d$ times the class of a point, where $d$ is (by
definition) the degree of $X$.
\ind  The following result is a direct
consequence of Property (1.3): \smallskip
 {\bf Proposition} 1$.-$ {\it
Let $X$ be a projective manifold, of dimension $n\ge 2$, of degree $d$.
Assume:} \indp (i) {\it The ordinary cohomology algebra $H^*(X,{\bf Q})$ is
spanned by $H$ and} $H^n(X,{\bf Q})$;\indp(ii) {\it One has $K_X=-kH$ with}
$k>{n\over 2}$. \indp(iii) {\it If $n=2k-1$, $H^n(X,{\bf Q})$ is nonzero.
\ind There exists an integer $\mu(X)$ such that the quantum cohomology
algebra\break $H^*(X,{\bf Q})$ is the algebra generated by $H$ and  $H^n(X,{\bf
Q})_{\rm o}$, with the  relations}: $$H^{n+1}=\mu(X)\,H^{n+1-k}\qquad
H\cdot\alpha=0\qquad \alpha\cdot\beta=(\alpha\,|\,\beta){1\over
d}\,(H^n-\mu(X)\,H^{n-k})\leqno(R)$$ {\it for} $\alpha,\beta\in H^n(X,{\bf
Q})_{\rm o}$.
\ind (Recall that the primitive cohomology $H^n(X,{\bf Q})_{\rm o}$ is
by definition equal to $H^n(X,{\bf Q})$ if $n$ is odd, and to the orthogonal of
$H_{n\over 2}$ if $n$ is even.)
  \ind Let
$p$ be an integer, with $|p|<k$. According to (1.3), one has $$H\cdot
H_{k+p-1}=H_{k+p}+\ell_p\,H_p\ ,\leqno(1.6)$$where $\ell _p$ is an integer
(which is zero for $p<0$).  Intersecting both sides with $H_{n-p}$ gives
$\ell_p={1\over d}\langle H,H_{n-p},H_{k+p-1}\rangle$ (so that $\ell
_p=\ell _{n-k+1-p}$). \ind From $(1.6)$ one obtains inductively, for $|p|<k$,
$$H_{k+p}=H^{k+p}-(\sum_{i=0}^p \ell _i)\,H^p\ .\leqno(1.7)$$If $n<2k-1$, we
can
apply this with $p=n-k+1$, yielding
$$0=H_{n+1}=H^{n+1}-(\sum_{i=0}^{n+1-k} \ell _i)\,H^{n+1-k}\ .$$ If
$n=2k-1$, one finds $\ H\cdot H_n=\ell _k\,H^k+m\ $ for some integer $m$.
If $m$ is nonzero $H$ is  invertible in $H^*(X,{\bf Q})$; since $H\cdot
H^n(X,{\bf Q})$ is zero for  degree reasons, this  implies $H^n(X,{\bf
Q})=0$.   Therefore under the hypotheses of the Proposition we obtain
 $$ H^{n+1}=\mu(X)\,H^{n+1-k}\qquad {\rm
with}\quad \mu(X)=\sum_{i=0}^{n+1-k} \ell _i\ .\leqno(1.8)$$
\ind Let $\alpha\in H^n(X,{\bf Q})_{\rm o}$. If $n\not=2k-2$ one has
again $H\cdot\alpha=0$ for degree reasons. If $n=2k-2$,
$H\cdot\alpha$ belongs to $H^0(X,{\bf Q})$ and one
has$$\eqalign{(H\cdot\alpha\,|\,H_n)&=(H\cdot H_n\,|\,\alpha)&\cr
&=((H_{n+1}+\ell _{n/2}\,H_{n/2})\,|\,\alpha)&\hbox{by (1.6)}\cr &=0
&\hbox{since $\alpha$ is primitive,} }$$hence again $H\cdot\alpha=0$.
\ind Let $\alpha,\beta\in H^n(X,{\bf Z})_{\rm o}$. By (1.3) and (1.7) there
exists an
integer $q$ such that $$\alpha\cdot\beta=(\alpha\,|\,\beta){1\over
d}H^n+q\,H^{n-k}$$
Multiplying by $H$ and using (1.8) yields
$\displaystyle q=-(\alpha\,|\,\beta){\mu(X)\over d}$ , which gives the
last relation (R).
\ind Finally we just have to remark  that the ${\bf Q}$\tx algebra
spanned by $H$ and $H^n(X,{\bf Q})_{\rm o}$ with the relations (R) has
 the same dimension as $H^*(X,{\bf Q})$, so that all relations follow
from (R). \cqfd
\bigskip
{\it Remarks}$.-$  1) Assume moreover that the variety of lines
contained in $X$ has the expected dimension $n+k-3$, and that $H$ is very
ample,
i.e. is the class of a hyperplane section of $X\i{\bf P}^N$. Then according to
(1.5)
  $d\,\ell _p$ is the number of lines in $X$ meeting two general
linear spaces of codimension $n-p$ and $k+p-1$ respectively. In
particular, $\ell _0$ is the number of lines passing through a point;
$\ell _1$ is the number of lines meeting a general line in $X$ (if $n\ge 3$).
\ind 2) Condition
(iii) is necessary. A counter-example is given by a general linear section of
codimension $3$ of the Grassmannian ${\bf G}(2,5)$. This is a Fano threefold of
index $k=2$, which satisfies the hypotheses (i) and (ii) of the Proposition.
For
such a threefold one has by (1.3)  $\ H\cdot H_3=\ell _0H_2+c$ , with
$c={1\over
d}\langle H,H_3,H_3\rangle^{}_2$. From $H^2=H_2+\ell _0\ $ and $\ H^3=H_3+(\ell
_0+\ell _1)H\ $ (1.7) we deduce $$H^4=(2\ell _0+\ell _1)\,H^2+c-\ell _0^2\ .$$
Easy
geometric computations give $\ell _0=3$, $\ell _1=5$, $c=10$, hence $c-\ell
_0^2=1\not=0$. \vskip1truecm
{\bf 2. Complete intersections}
\smallskip
\ind Let $X$ be a smooth complete intersection in ${\bf P}^{n+r}$ of degree
$(d_1,\ldots,d_r)$ and dimension $n\ge 2$, with
 $n\ge 2\sum (d_i-1)-1$. One has $K_X=-kH$,
with $k=$ $n+1-\sum (d_i-1)$; therefore the inequality on $n$  ensures that
condition (ii) of Proposition 1  holds. The space $H^n(X,{\bf Q})$ is  nonzero
except for odd-dimensional quadrics [D], so condition (iii) holds as well;
finally
(i) holds by the weak Lefschetz theorem. Therefore the quantum cohomology of
$X$ is
given by Proposition 1; to achieve the proof of the Theorem it remains to
compute the
number $\mu(X)=\sum\ell _p$.
 \ind In view of $(1.1)$ we can assume that $X$ is general; then the variety of
lines
(resp.\ conics, resp.\ twisted cubics) contained in $X$ has the expected
dimension:
see for instance [E-S], where the proof (given  for the case of twisted cubics)
adapts immediately to the easier cases of lines and conics.  Therefore $d\ell
_p$ is
the number of lines in $X$ meeting  two general linear spaces of codimension
$n-p$
and $k+p-1$ respectively (Remark 1). \ind Let $V$ be a complex vector space, of
dimension $N$; let us denote by $G={\bf G}(2,V)$ the Grassmannian of lines
in the projective space ${\bf P}(V)${\footnote{\parindent
0.4cm$^1$}{\vtop{\hsize15.5truecm\noindent \eightpoint
We use the naive convention, i.e. ${\bf P}(V)$ is the variety of lines
in $V$.}}\parindent 0cm}. On $G$ we have a tautological exact sequence
$$0\rightarrow
S\longrightarrow {\cal O}_G\otimes_{\bf C}V\longrightarrow Q\rightarrow 0\ ,$$
where  the   sub- and quotient bundles $S$ and $Q$ are of rank $2$ and
$N-2$ respectively. \ind The Chern classes $c_1,\ldots,c_{N-2}$ of $Q$ are
represented by the {\it special Schubert cycles}:  $$c_p={\rm cl}\,\{\ell \in
G\ |\
\ell \cap H_{p+1}\not=\emptyset\}$$ where $H_{p+1}$ is a fixed linear subspace
of
${\bf P}(V)$ of codimension $p+1$. In particular, the subvariety of lines in
$G$
meeting two general linear spaces of codimension $p+1$ and $q+1$ has cohomology
class
$c_p\,c_q$.  \ind Let $f\in \sy^dV^*$ be a homogeneous polynomial of degree $d$
on
${\bf P}(V)$. It defines by restriction a global section $\overline{f}$ of
$\sy^dS^*$, which vanishes exactly at the points of $G$ where the corresponding
line  is contained in the hypersurface $f=0$. In other words, the subvariety of
lines contained in this hypersurface is the zero locus of $\overline{f}\in
H^0(G,\sy^dS^*)$. If $f$ is general enough, it has the expected codimension
$d+1$,
and therefore its cohomology class is the top Chern class  $c_{d+1}(\sy^dS^*)$.
Hence the cohomology class of the variety  of lines contained in our complete
intersection $X$ is
 $c_{d_1+1}(\sy^{d_1}S^*)\ldots c_{d_r+1}(\sy^{d_r}S^*)$.
 Therefore we
find $$\ell _p={1\over d}\int_Gc_{d_1+1}(\sy^{d_1}S^*)\ldots
c_{d_r+1}(\sy^{d_r}S^*)\, c_{n-1-p}\,c_{k-2+p}$$
(recall that $k=n+r+1-\sum d_i$).
\ind We will compute this number using the Chern classes
$x=c_1(S^*)$, $y=c_2(S^*)$, or rather the virtual classes
$\alpha,\beta$ such that $x=\alpha+\beta$, $y=\alpha\beta$. The
Schubert cycles $c_p$ are then given by
$$\eqalign{
1+c_1+\ldots+c_{N-2}=(1-x+y)^{-1}&=(1-\alpha)^{-1}(1-\beta)^{-1}\cr
&={1\over \alpha-\beta}\ \Bigl({\alpha\over 1-\alpha}-{\beta\over
1-\beta}\Bigr)\ ,}$$ hence $\displaystyle c_p={\alpha^{p+1}-\beta^{p+1}\over
\alpha-\beta}$; the Chern class $c_{d+1}(\sy^{d}S^*)$ is equal to
$\displaystyle \prod_{j=0}^d(j\alpha+(d-j)\beta)$. To integrate we use the
following lemma: \smallskip  {\bf Lemma}$.-$ {\it Let $P\in {\bf
C}[\alpha,\beta]$ be a symmetric
homogeneous polynomial of degree $2(N-2)$ {\rm (}so that
$P(\alpha,\beta)$ is a polynomial of maximum degree in the Chern classes
$x$ and $y)$. Then $\int_GP(\alpha,\beta)$ is the
coefficient of $\alpha^{N-1}\beta^{N-1}$ in
$-{1\over 2}(\alpha-\beta)^2\,P(\alpha,\beta)$.}

\ind This is probably well-known; let me give a quick proof for the sake of
completeness. Put $\displaystyle c_p={\alpha^{p+1}-\beta^{p+1}\over
\alpha-\beta}$ for {\it all} $p$. The (usual!) cohomology  algebra of $G$ is
the algebra of symmetric polynomials in $\alpha,\beta$, modulo
the ideal generated by $c_{N-1}$ and $c_N$ [G].
Consider the linear form which associates to a symmetric polynomial
$P(\alpha,\beta)$  the coefficient of $\alpha^{N-1}\beta^{N-1}$ in
$-{1\over 2}(\alpha-\beta)^2\,P(\alpha,\beta)$. It vanishes on the ideal
$(c_{N-1},c_N)$ and on the polynomials of degree $<2N-4$, hence  factors
through a linear form $\ell :H^{2N-4}(G,{\bf Q})\rightarrow {\bf Q}$,
necessarily proportional to $\int_G$. Let us evaluate these two forms on
the polynomial $c^2_{N-2}$. One has
$(\alpha-\beta)^2\,c^2_{N-2}=(\alpha^{N-1}-\beta^{N-1})^2$, hence $\ell
(c^2_{N-2})=1$; on the other hand $\int_Gc^2_{N-2}$ is the number of
lines in ${\bf P}(V)$ through 2 points, that is $1$. This proves the lemma.
\medskip
\ind  Let us apply the lemma to the polynomial
$\ F(\alpha,\beta) \,c_{n-1-p}\,c_{k-2+p}$ , where $\displaystyle
F(\alpha,\beta)=\sum_{j=1}^{\delta-1} a_j\,\alpha^j\beta^{\delta-j}$ is a
symmetric homogeneous polynomial of degree $\delta:=$ $\sum(d_i+1)$. One
has  $$\eqalign{(\alpha-\beta)^2\,c_{n-1-p}\,c_{k-2+p}&=(\alpha^{n-p}-
\beta^{n-p})\,(\alpha^{k-1+p}-\beta^{k-1+p})\cr
&=\alpha^{n+k-1}+\beta^{n+k-1} -\alpha^{n-p}\beta^{k-1+p}
-\alpha^{k-1+p}\beta^{n-p}}\ .$$
Since
$N=n+r+1$, the coefficient of $\alpha^{N-1} \beta^{N-1}$ in
$(\alpha-\beta)^2\,F(\alpha,\beta)\,c_{n-1-p}\,c_{k-2+p}$ is
$2a_{r-k+1}-2a_{r+p}$;
if moreover  $F(\alpha,\beta)$ is divisible by
$(\alpha\beta)^r$, the first coefficient is zero (recall that
$k>{n\over 2}\ge 1$). Applying this to the polynomial $F(\alpha,\beta)=$ $
c_{d_1+1}(\sy^{d_1}S^*)\ldots
c_{d_r+1}(\sy^{d_r}S^*)$ we get $\ell _p=a_{r+p}$, that is
$$\sum_{p=0}^{n+1-k}\ell _p \,\alpha^{r+p}\beta^{\delta-r-p}
={1\over d}\,\prod_{i=1}^r\,\prod_{j=0}^{d_i}(j\alpha+(d_i-j)\beta)\
.\leqno(2.1)$$
Taking $\alpha=\beta=1$ gives $\displaystyle \mu(X)=\sum\ell
_p=\prod_{i=1}^rd_i^{d_i}$, which achieves the proof of the Theorem. Note
that formula $(2.1)$ gives explicit expressions for the $\ell _p$'s, for
instance
$$\ell _0=\prod_{i=1}^rd_i!\leqno(2.2)$$
$$\ell _1=\prod_{i=1}^rd_i!\ \Bigl(\sum_{1\le i\le r\atop 1\le
j<d_i}{d_i-j\over j}\Bigr)\ ,\leqno(2.3)$$and so on.
 \vskip1truecm

{\bf 3. Application I: enumerative formulas}
\smallskip
\ind Let $X$ be a smooth projective manifold satifying the hypotheses of
Proposition
1; it follows from that Proposition that all the triple products  $\langle
H_p,H_q,H_r\rangle^{}_i$ can be computed in terms of the integers $\ell
_p$. If the variety of lines, conics or twisted cubics in $X$ has the
expected dimension, this gives some nice enumerative formulas which we are
going to describe.
\ind Let $p,q,r$ be
positive integers $\le n$ such that $p+q+r=n+k$; we arrange them so that
$p\le q\le r$. Since $2k>n$ by hypothesis this implies $p<k$ and
$k\le p+q<2k$ (hence  $q\ge k$). Therefore $$H_p\cdot
H_q=H^p\cdot(H^q-(\sum_{i=0}^{q-k}\ell _i)\,H^{q-k})=
(\mu(X)-\sum_{i=0}^{q-k}\ell _i)\ H_{p+q-k}\ ,$$
$$\langle H_p,H_q,H_r\rangle^{}_1=(H_p\cdot H_q\,|\,H_r)=d\,
(\mu(X)-\sum_{i=0}^{q-k}\ell _i)\ .\leqno{\rm
hence}$$
Using the
equalities $\ell _i=\ell _{n+1-k-i}$ and the convention
$\ell _i=0$ for $i>n+1-k$, we find\smallskip

{\bf Proposition} 2$.-$ {\it Assume that the variety of lines
contained in $X$ has the expected dimension $n+k-3$. Let $p,q,r$ be
positive integers such that $p\le q\le r\le n$ and $p+q+r=n+2k$. The number
of lines in $X$ meeting three general linear spaces of codimension $p$,$q$
and $r$ respectively in ${\bf P}^{n+r}$ is $\displaystyle
d\,\sum_{i=0}^{n-q}\ell _i\ .$}
\ind Actually this could also be obtained by a computation in the
Grassmannian as in \S 2. This is  probably also the case for the next
results, though the computation would  be much more involved.
 \ind Let us look at conics.
Let $p,q,r$ be positive integers such that $p+q+r=$ $n+2k$; as above we
assume $p\le q\le r\le n$. Moreover we will assume $k<n$,
which  excludes only the trivial case of quadrics  [K-O].
This implies $p<k$ and therefore  $2k<p+q<3k$. We have as before
$\displaystyle\ H_p\cdot H_q=(\sum_{i=0}^{n-q}\ell _i)\ H^{p+q-k}$ ; since
$ H^{p+q-k}=$ $\displaystyle H_{p+q-k}+(\sum_{j=0}^{n-r}\ell _j)\,
H_{n-r}$, we obtain
\vskip-5pt
$$\langle H_p,H_q,H_r\rangle^{}_2=(H_p\cdot
H_q\,|\,H_r)=d\, (\sum_{i=0}^{n-q}\ell _i)\ (\sum_{j=0}^{n-r}\ell
_j)\ .$$
{\bf Proposition} 3$.-$ {\it Assume that $X$ is not a quadric, and that
the variety of conics contained in $X$ has the expected dimension $n+2k-3$.
Let $p,q,r$ be positive integers such that $p\le q\le r\le n$ and
$p+q+r=n+2k$. The number of conics in $X$ meeting three general linear
spaces of codimension $p$,$q$ and $r$ respectively in ${\bf P}^{n+r}$ is
$\displaystyle d\,(\sum_{i=0}^{n-q}\ell _i)\,(\sum_{j=0}^{n-r}\ell _i)\ .$}
\smallskip \ind This has to be taken with a grain of salt in the case $p=1$,
$q=r=n$, because every hyperplane meets a conic twice, so the above number must
be
divided by $2$. Since $H_n$ is $d$ times the class of a point we find that the
number of conics in $X$  through $2$ general points is ${\ell _0^2\over 2d}$,
where
$\ell _0$ is the number of lines through a general point in the intersection of
$X$
with a general linear space of codimension $k-2$. For complete intersections
formula
(2.2) gives: \smallskip {\bf Corollary}$.-$ {\it Let $X$ be a
smooth complete intersection of degree $(d_1,\ldots,d_r)$ in ${\bf
P}^{n+r}$, with  $n= 2\sum
(d_i-1)-1$.  The number of conics in $X$
passing through $2$ general points is $\displaystyle {1\over
2d}\prod_{i=1}^r(d_i!)^2$.} \medskip {\bf Example}$.-$ Let $X$ be a cubic
threefold, $P,Q$  two general
points in $X$, $L,M$ two general lines. We find that there are $6$
conics in $X$ through $P$ and $Q$ -- a fact that can easily be checked
geometrically (the line $\langle P,Q\rangle$ meets $X$ along a third point
$R$; conics in $X$ through $P$ and $Q$ are in one-to-one correspondence
with lines through $R$). Similarly from Proposition 3 we find $14$ conics
through $P$ meeting $L$ and $M$.\medskip

 \ind The
computation for twisted cubics is very similar. Let $p,q,r$ be
positive integers  with $p\le q\le r\le n$ and $p+q+r=n+3k$. Since $2k>n$
this implies $p\ge k$ and $p+q>k+n$. We have $$\eqalign{H_p\cdot
H_q&=(H^p-(\sum_{i=0}^{p-k}\ell
_i)\,H^{p-k})\cdot(H^q-(\sum_{j=0}^{q-k}\ell _j)\,H^{q-k})\cr
&=\bigl(\mu(X)^2-\mu(X)\sum_{j=0}^{q-k}\ell _j-\mu(X)
\sum_{i=0}^{p-k}\ell_i+(\sum_{i=0}^{p-k}\ell
_i)(\sum_{j=0}^{q-k}\ell_j)\bigr)H^{p+q-2k}\cr
&=\bigl(\mu(X)-\sum_{i=0}^{p-k}\ell
_i)(\mu(X)-\sum_{j=0}^{q-k}\ell
_j)\bigr)\ \bigl(H_{p+q-2k}+(\sum_{m=0}^{p+q-3k}\ell
_m)\,H_{p+q-3k}\bigr)
}$$
Reasoning as above we get:\smallskip
{\bf Proposition} 4$.-$ {\it  Assume that the variety of twisted cubics
contained in $X$ has the expected dimension $n+3k-3$. Let $p,q,r$ be
positive integers such that $p\le q\le r$ $\le n$ and $p+q+r=n+3k$. The number
of twisted cubics
 in $X$ meeting three general linear spaces of codimension $p$,$q$
and $r$ respectively in ${\bf P}^{n+r}$ is $$\displaystyle
d\
(\sum_{i=0}^{n-p}\ell _i)\ (\sum_{j=0}^{n-q}\ell
_j)\ (\sum_{m=0}^{n-r}\ell_m)\ .$$}
\ind  In particular:
\smallskip
{\bf Corollary}$.-$ {\it Let $X$ be a
smooth complete intersection of degree $(d_1,\ldots,d_r)$ in ${\bf
P}^{n+r}$, with $n=3\sum(d_i-1)-3$. Then
the number of twisted cubics in $X$ passing through $3$ general points
is $\displaystyle {1\over d^2}\prod(d_i!)^3$.}
\medskip
{\bf Example}$.-$ Going back to our cubic threefold we find that the
number of twisted cubics through $3$
general points is $24$. I do not know an easy geometric way to get this
number (G. Ellingsrud checked it using the method of [E-S] and a
computer).
\vskip1truecm
{\bf 4. Application II: the primitive cohomology}
\smallskip
\ind So far we have only considered the subalgebra of $H^*(X,{\bf Z})$
generated by $H$. In this last section I would like to look at the remaining
part. Because of the relations $(R)$, the only interesting triple product
which appears is  $\langle H_k,\alpha,\beta\rangle^{}_1$ for $\alpha,\beta\in
H^n(X,{\bf Z})_{\rm o}$.  Since $H_k\cdot\alpha=(H^k-\ell _0)\cdot\alpha=-\ell
_0\alpha$, we get $$\langle
H_k,\alpha,\beta\rangle=(H_k\cdot\alpha\,|\,\beta)=-\ell _0
\,(\alpha\,|\,\beta)\ .\leqno(4.1)$$ \ind Let us translate this geometrically
using (1.4). We suppose given a smooth subvariety $Y$ of $X$, of codimension
$k$
and degree $d_Y$, such that  the variety  $\Gamma$ of lines in $X$ meeting $Y$
is smooth, of dimension $n-2$. For instance we can take for $Y$ a general
linear section of codimension $k$ in $X$; if $k=n-1$, we can take for $Y$ a
line. The correspondence
\vskip-10pt$$\hbox{$\diagram{R&\phfl{q}{}&X\cr \vfl{p}{}&&\cr
\Gamma&&\cr}$\hskip2cm with $\ R:=\{(L,x)\in\Gamma\times X\ |\ x\in
L\}$}$$\vskip-10pt
gives rise to a homomorphism $\varphi=p_*q^*:H^n(X,{\bf Z})\rightarrow
H^{n-2}(\Gamma,{\bf Z})$. By definition this is a morphism of Hodge
structures, i.e.  $\varphi^{}_{\bf C}$ maps $H^{p,q}(X)$ into
$H^{p-1,q-1}(\Gamma)$ for $p+q=n$.   \medskip
  {\bf Proposition} 5$.-$ {\it One has $\displaystyle
(\varphi(\alpha)\,|\,\varphi(\beta))=-{\ell _0d_Y\over d}\,(\alpha\,|\,\beta)$
for $\alpha,\beta\in H^n(X,{\bf Z})_{\rm o}$}.

\ind The moduli space ${\cal M}_1$ of degree $1$ maps ${\bf
P}^1\rightarrow X$ has a  natural smooth compactification, namely
$\overline{\cal M}_1=U\times_FU\times_FU$, where   $F$ is the
variety of lines in $X$ and  $U$ the tautological family of lines over $F$; the
map $e_i:\overline{\cal M}_1\rightarrow X\quad(0\le i\le 2)$ is obtained by
composing
 the projection $p_{i+1}$ with  the natural map $U\rightarrow X$.
The inverse
image of  $Y$ under  $e_0$ is then identified with the
fibered product $R\times_\Gamma R$, in such a way that the evaluation map
$e_i:R\times_\Gamma R\rightarrow X$ is $q\rond p_i$. (1.4) yields
$$\langle
Y,\alpha,\beta\rangle^{}_1=\int_{R\times_\Gamma R}e_1^*\alpha\  e_2^*\beta\ .$$
\ind   Since $R$ is a ${\bf P}^1$\tx bundle over $\Gamma$ and the class $q^*H$
is
transversal to the fibres, the map $\lambda:H^{n-2}(\Gamma,{\bf Z})\oplus
H^n(\Gamma,{\bf Z})\longrightarrow  H^n(R,{\bf Z})$ given by
$\lambda(\gamma,\delta)= p^*\gamma\,.\, q^*H+p^*\delta$ is an isomorphism,
which
satisfies $p_*\lambda(\gamma,\delta)=\gamma$. Let us write
$$q^*\alpha=p^*\varphi(\alpha)\,.\,q^*H+p^*\alpha'\qquad,\qquad
q^*\beta=p^*\varphi(\beta)\,.\,q^*H+p^*\beta'\ .$$ Let $\pi=p\rond
p_1=p\rond p_2$ be the projection of $R\times_\Gamma R$ onto $\Gamma$. One
has
$$e_i^*\alpha=p_i^*q^*\alpha=\pi^*\varphi(\alpha)\,p_i^*q^*H+\pi^*\alpha'
\quad\hbox{ , and similarly for }\ e_2^*\beta\ .$$ For degree reasons
the last terms disappear in the product $e_1^*\alpha\,.\,e_2^*\beta$, and we
get

$$\displaystyle \langle
Y,\alpha,\beta\rangle^{}_1=(\varphi(\alpha)\,|\,\varphi(\beta))\int_{L\times
L}p_1^*q^*H\,.\,p_2^*q^*H\ ,$$where $L$ is a general line
intersecting $Y$. The value of the integral is obviously
 $1$; since the cohomology class of $Y$ is $\displaystyle {d_Y\over d}H_k$,
the result follows from $(4.1)$.
\medskip
{\bf Example}$.-$ Let us go back to our favorite example, the cubic
threefold, taking for $Y$ a generic line in $X$.  Then $\Gamma$ is a smooth
curve; the map $\varphi:H^3(X,{\bf Z})\rightarrow H^1(\Gamma,{\bf Z})$
gives rise to a morphism $\Phi:JX\rightarrow J\Gamma$, where $J\Gamma$ is the
Jacobian of $\Gamma$ and $JX$ the {\it intermediate Jacobian} of $X$ (see
e.g. [C-G]);   the formula
$(\varphi(\alpha)\,|\,\varphi(\beta))= -2\,(\alpha\,|\,\beta)$ for
$\alpha,\beta\in H^3(X,{\bf Z})$ given by Proposition 5 means  that {\it the
principal polarization of $J\Gamma$ induces twice the principal polarization
of $JX$}. One deduces easily from this that {\it the intermediate Jacobian
$JX$ is isomorphic} (as a principally polarized Abelian variety) {\it to the
Prym variety associated to $\Gamma$ and the natural involution
 of} $\Gamma$ which
maps a line $L$ to the third line cut down on $X$ by the $2$\tx plane
spanned by $Y$ and $L$ -- a fundamental fact for the geometry of the cubic
threefold, due to Mumford (see Appendix C of [C-G]).

\baselineskip14pt
\vskip 2cm
\def\num#1{\item{\hbox to\parindent{\enskip [#1]\hfill}}}
\def\pc#1{\tenrm#1\sevenrm}
\parindent=1.5cm
\centerline{\bf REFERENCES}
\bigskip
\num{C-G} 	H. {\pc CLEMENS}, P. {\pc GRIFFITHS}: {\sl The intermediate
Jacobian of the cubic threefold}. Ann. of Math. {\bf 95}, 281-356 (1972).
\smallskip
\num{C-M} B. {\pc CRAUDER}, R. {\pc MIRANDA}: {\sl Quantum cohomology of
rational
surfaces}. Preprint alg-geom/9410028 (1994). \smallskip
\num{D} P. {\pc DELIGNE}: {\sl Cohomologie des intersections compl\`etes}. SGA
7,
Expos\'e XI;  Lecture Notes in Math. {\bf 340}, Springer-Verlag (1973).
\smallskip \num{E-S} G. {\pc ELLINGSRUD}, S.-A. {\pc STROMME}: {\sl Bott's
formula
and enumerative geometry}. Preprint alg-geom/9411005 (1994). \smallskip
\num{G} A. {\pc GROTHENDIECK}: {\sl Sur quelques propri\'et\'es
fondamentales en th\'eorie des intersections}. S\'eminaire Chevalley 1958
(``Anneaux
de Chow et applications"), exp. IV. IHP (1958). \smallskip \num{K-O} S. {\pc
KOBAYASHI}, T. {\pc OCHIAI}: {\sl Characterizations of complex projective
spaces and
hyperquadrics}. J. Math. Kyoto Univ. {\bf 13}, 31-47 (1973). \smallskip
\num{R-T} Y.
{\pc RUAN}, G. {\pc TIAN}: {\sl A mathematical theory of quantum cohomology}.
Pre\-print (1995). \smallskip \num{V} C. {\pc VAFA}: {\sl Topological mirrors
and
quantum rings}, Essays on mirror mani\-folds (S.-T. Yau ed.), 96-119;
International
Press, Hong-Kong (1992). \smallskip
\num{W} E. {\pc WITTEN}: {\sl Topological sigma models}. Commun. Math. Phys.
{\bf
118}, 411-449 (1988).
\vskip1cm
\hfill\hbox to 5cm{\hfill A. Beauville\hfill}\par
\hfill\hbox to 5cm{\hfill URA 752 du CNRS\hfill}\par
\hfill\hbox to 5cm{\hfill Math\'ematiques -- B\^at. 425\hfill}\par
\hfill\hbox to 5cm{\hfill Universit\'e Paris-Sud\hfill}\par
\hfill\hbox to 5cm{\hfill 91 405 {\pc ORSAY} Cedex, France\hfill}\par
 \bye